\documentclass[12pt,reqno]{article}
\usepackage{amsmath,amsthm,amssymb}

\numberwithin{equation}{section}

\usepackage{tikz-cd}
\usepackage[all]{xy}
\usepackage{mathtext}
\usepackage[T1,T2A]{fontenc}
\usepackage[utf8]{inputenc}
\usepackage[russian,english]{babel}
\usepackage{amsfonts,longtable,amssymb,amsmath,latexsym, dsfont}
\usepackage{wrapfig}
\usepackage{xypic}
\usepackage{mathtools}
\usepackage{hyperref}
\title{Selfsimilar Hessian and conformally K\"ahler manifolds}
\author{Pavel Osipov\footnote{National Research University Higher School of Economics, Russian Federation }\footnote{Pavel Osipov is 
		partially supported by the HSE University Basic
		Research Program, Simons Foundation, and by the contest “Young Russian
		Mathematics”.} 
}

\usepackage{amsthm, amsmath, amssymb, amsfonts, amscd, amscd, geometry}
\newsavebox{\ssa}
\renewcommand{\d}{\partial}

\newcommand{\dxi}{\frac{\partial}{\partial x^i}}

\newcommand{\R}{\mathbb R}

\newcommand{\CC}{\mathbb{C}}

\renewcommand{\L}{\mathcal{L}}

\theoremstyle{definition}
\newtheorem{theorem}{Theorem}[section]
\newtheorem{lemma}[theorem]{Lemma}
\newtheorem{proposition}[theorem]{Proposition}

\newtheorem{cor}[theorem]{Corollary}
\newtheorem{defin}[theorem]{Definition}

\theoremstyle{remark}
\newtheorem{rem}{Remark}[section]
\usepackage{indentfirst}

\begin{document}
	
		\maketitle
		
		\begin{abstract}
		Let $(M,\nabla,g)$ be a Hessian manifold. Then the total space of the tangent bundle $TM$ can be endowed with a K\"ahler structure $\left(I,{g^\text{r}}\right)$. We say that a homogeneous Hessian manifold is a Hessian manifold $(M,\nabla,g)$ endowed with a transitive action of a group $G$ preserving $\nabla$ and $g$.   We construct by a Hessian (special K\"ahler) structure on a simply connected manifold with a certain condition a K\"ahler (hyeper-K\"ahler) structure on the tangent (cotangent) bundle. A selfsimilar Hessian (K\"ahler) manifold is a Hessian manifold endowed with a homothetic vector field $\xi$. We construct by a selfsimilar Hessian (K\"ahler) structure on a simply connected manifold with a certain condition a conformally K\"ahler (hyeper-K\"ahler) structure on the tangent (cotangent) bundle.
	\end{abstract}

		\tableofcontents

	\section{Introduction}
	An open cone $V\subset \R^n$ is called {\bfseries regular} if it does not contain any straight full line. The tube neighborhood $V+\sqrt{-1} \R^n\subset \CC^n$ is biholomorphic to a bounded complex domain. All complex domains arising by this way are called {\bfseries Siegel domains of the first kind}. All complex affine automorphisms of $V+\sqrt{-1} \R^n\subset \CC^n$ can be  written in  the  form $x+\sqrt{-1} y \to Ax +\sqrt {-1}\left(Ax+b\right)$ where $A$ is a linear automorphism of $V$ and $b\in \R^n$. We consider linear automorphisms of cones and complex affine automorphisms of Siegel domains of the first kind. In these assumptions, $V+\sqrt{-1}\R^n$ is homogeneous if and only if $V$ is homogeneous. Any homogeneous Siegel domain of the first kind admits an invariant K\"ahler structure (see \cite{VGP} or \cite{C}).
	
	In this paper we modify the construction of invariant K\"ahler structures on Siegel domains of the first kind and get a construction of a certain class of homogeneous conformally K\"ahler manifolds. In particular, we confirm that homogeneous Siegel domains of the first kind admit invariant conformally K\"ahler structures.
	
	A {\bfseries flat affine manifold} is a differentiable manifold equipped with a flat, torsion-free connection. Equivalently, it is a manifold equipped with an atlas such that all translation maps between charts are affine transformations (see \cite{FGH} or \cite{shima}). 
	A {\bfseries Hessian manifold} is an affine manifold with a Riemannian metric which is locally equivalent to a Hessian of a function. Any Kähler metric can be locally defined as a complex Hessian $\d \bar \d \varphi$. Thus, the Hessian geometry is a real analogue of the Kähler one.

	Hessian manifolds have many different applications: in supersymmetry (\cite{CMMS}, \cite{CM}, \cite{AC}), in convex programming
	(\cite{N}, \cite{NN}), in the Monge-Ampère Equation (\cite{F1}, \cite{F2}, \cite{G}), in the WDVV equations (\cite{T}).

	A Kähler structure $(I,{g^\text{r}})$ on $TM$ can be constructed by a Hessian structure $(\nabla,g)$ on $M$ (see \cite{shima}). The correspondence 
	$$
	\text{r}:\{\text{Hessian manifolds}\} \to \{\text{K\"ahler manifolds}\}
	$$
	$$
\ \ \ \	\ (M,\nabla,g)\  \ \to \ \ (TM,I,{g^\text{r}})
	$$
	is called the {\bfseries (affine) r-map}. In particular, this map associates some special Kähler manifolds to special real manifolds (see \cite {AC}). In this case, r-map describes a correspondence between the scalar geometries for supersymmetric theories in dimension {5 and 4.} See \cite{CMMS} for details on the r-map and supersymmetry.
	
	Any regular cone admits a function $\varphi$ called {\bfseries characteristic function} such that $g_{can}=\text{Hess} \left(\ln \varphi\right)$ is a Hessian metric which is  invariant with respect to all automorphisms of the cone. (\cite{vinb}). The r-map constructs an invariant K\"ahler structure $\left(I, {g^\text{r}_{can}}\right)$ on $TV\simeq V\oplus \sqrt{-1} \R^n$. Thus, any homogeneous Siegel domain of the first kind admits an invariant K\"ahler structure. The K\"ahler potential of ${g^\text{r}_{can}}$ equals $4\pi^*\left(\ln \varphi\right)$.

	The construction of the invariant K\"ahler structure on $V\oplus \sqrt{-1}\R^n$ is well known (see \cite{VGP} or \cite{C}). We modify this construction. A {\bfseries (globally) conformally K\"ahler manifold} $\left(M,I,\omega\right)$ is a complex manifold endowed with Riemannian metric $g$ which is (globally) conformally equivalent to a Kahler one. We consider the metric $g_{con}=\text{Hess}\ \varphi$ on a regular homogeneous cone $V$. This metric is invariant under $\text{Aut}(V)\cap \text{SL}(\R^n)$ and coincides with $g_{can}$ on the hypersurface $\{\varphi(x)=1\}$. The dilation $x \mapsto qx$ acts on  $g_{con}$ by $\lambda_q^* g_{con} = q^{-n} g_{con}$. The K\"ahler metric ${g^\text{r}}_{con}$ on $V\oplus \sqrt {-1} \R^n$ constructed by the r-map is not invariant but it is conformally equivalent to the invariant Riemannian metric $r^{-2}  {g^\text{r}}_{con}$ on the homogeneous domain $V\times \sqrt{-1} \R^n$. Thus, Siegel domains of the first kind admit two different invariant structures: K\"ahler and conformally K\"ahler.

	We generalize this construction  to  selfsimilar  Hessian manifolds. A {\bfseries selfsimilar Hessian manifold} is a Hessian manifold endowed with a homothetic vector field generating a flow of affine automorphisms. For example, an $n$-dimansional regular convex cone $\left(V,g_{con}\right)$ endowed with a field $-\frac{2}{n}\sum x^i\dxi$ is a selfsimilar Hessian manifold. 
	
	The main result of the paper is the following. 
	\begin{theorem}
	Let $(M,\nabla,g,\xi)$ be a simply connected selfsimilar Hessian manifold such that $\xi$ is complete and $G$ be the group of affine isometries of $(M,\nabla,g)$ preserving $\xi$. Suppose that $G$ acts simply transitively on the level line $\{g(\xi,\xi)=1\}$. Then $TM$ admits a homogeneous conformally K\"ahler structure.  
\end{theorem} 
	Also, we adapt our construction for the c-map. 
		A  {\bfseries special K\"ahler manifold} $\left(M,I,g,\nabla\right)$ is a K\"ahler manifold  $\left(M,I,g\right)$ endowed with a torsion free symplectic connection $\nabla$ such that such that $g$ is a Hessian metric with respect to $g$ (see \cite{ACD}).

	A hyper-Kähler structure $\left(I_1,I_2,I_3, {g^\text{c}}\right)$ on $T^*M$ can be constructed by a special K\"ahler structure $\left(I,g,\nabla\right)$ on $M$. The correspondence 
	$$
	\text{c}:\{\text{Special K\"ahler manifolds manifolds}\} \to \{\text{hyper-K\"ahler manifolds}\}
	$$
	$$
	\ \ \ \	\ \left(M,I,g,\nabla\right)\  \ \to \ \ \left(T^*M,I_1,I_2,I_3, {g^\text{c}}\right)
	$$
	is called the {\bfseries (affine) c-map}. r-map describes a correspondence between the scalar geometries for supersymmetric theories in dimension {4 and 3}. See \cite{CMMS2} and \cite{ACM}  for details on the r-map
	
	A {\bfseries selfsimilar special K\"ahler manifold} is a special K\"ahler manifold $\left(M,I,g,\nabla\right)$  endowed with an affine vector field $\xi$ satisfying $\L_\xi g= 2g$. For example, any conic special K\"ahler manifold is selfsimilar special K\"ahler (see \cite{ACM} for the definition).

	  \begin{theorem}
	  	Let $\left(M,I,g,\nabla,\xi\right)$ be a selfsimilar special K\"ahler manifold, such that $\xi$ is complete and $G$ the group of affine holomorphic isometries of  $\left(M,I,g,\nabla\right)$ preserving $\xi$. Suppose that $G$ acts simply transitively on the level line $\{g(\xi,\xi)=1\}$. Then $T^*M$ admits a homogeneous conformally hyper-K\"ahler structure. 
	  	
	  \end{theorem}

	\section{Hessian and Kähler structures}
	\begin{defin}
		 A {\bfseries flat affine manifold} is a differentiable manifold equipped with a flat, torsion-free connection. Equivalently, it is a manifold equipped with an atlas such that all translation maps between charts are affine transformations (see \cite{FGH}).  
	\end{defin}

	\begin{defin}
		A Riemannianian metric $g$ on a flat affine manifold $(M,\nabla)$ is called to be a {\bfseries Hessian metric} if $g$ is locally expressed by a Hessian of a function
		$$
		g=\text{Hess} \ \varphi =\frac{\partial^2\varphi}{\partial x^i \partial x^j} dx^i dx^j,
		$$
		where $x^1,\ldots, x^n$  are flat local coordinates. A {\bfseries Hessian manifold} $(M,\nabla,g))$ is a flat affine manifold $(M,\nabla)$ endowed with a Hessian metric $g$. (see \cite{shima}). 
	\end{defin}
	
	Let $U$ be an open chart on a flat affine manifold $M$, functions $x^1,\ldots, x^n$ be affine coordinates on $U$, and $x^1,\ldots, x^n, y^1, \ldots, y^n$ be the corresponding coordinates on $TU$. Define the complex structure $I$ by $I(\frac{\partial}{\partial x^i})=\frac {\partial} {\partial y^i}$. Corresponding complex coordinates are given by $z^i=x^i+\sqrt {-1}y^i$. The complex structure $I$ does not depend on a choice of flat coordinates on $U$. Thus, in this way, we get a complex structure on the $TM$. 
	
	Let $\pi : TM \to M$ be a natural projection. Consider a Riemannian metric $g$ on $M$ given locally by
	\begin{equation*}\label{1}
	g_{i,j} dx^idx^j.
	\end{equation*}
	
	Define a bilinear form ${g^\text{r}}$ on $TM$ by
	\begin{equation*}\label{2}
	{g^\text{r}}=\pi^* g_{i,j} \left(dx^idx^j+dy^idy^j\right)
	\end{equation*}
	or, equivalently,
	\begin{equation}\label{101}
	{g^\text{r}}(X,Y)=\left(\pi^*g\right)(X,Y)+\left(\pi^*g\right)(IX,IY),
	\end{equation}
	for any $X,Y\in T\left(TM\right)$.
	
	\begin{proposition}[\cite{shima}, \cite{AC}] \label{2.3}
		Let $M$ be a flat affine manifold, $g$ and ${g^\text{r}}$ as above. Then the following conditions are equivalent:
			\begin{itemize}
				\item[(i)]$g$ is a Hessian metric.
		
		\item[(ii)] ${g^\text{r}}$ is a Kähler metric. 
	\end{itemize}
		Moreover, if
		$
		g=\text{Hess} \varphi
		$ 
		locally then ${g^\text{r}}$ is equal to a complex Hessian
		$$
		{g^\text{r}}=\text{Hess}_\mathbb{C} (4\pi^*\varphi)=\d\bar{\d}(4\pi^*\varphi).
		$$
	\end{proposition}

	\begin{defin}
		The metric ${g^\text{r}}$ is called a {\bfseries Kähler metric associated to $g$}. The correspondence which associates the Kähler manifold $(TM,{g^\text{r}})$ to a Hessian manifold $(M, g)$ is called the {\bfseries (affine) r-map} (see \cite{AC}).
	\end{defin}

	\section{The r-map for homogenous manifolds}

	\begin{defin}\label{31}
		Let $G$ be a group of affine automorphisms of a simply connected affine manifold $M$ and $\phi$ be the corresponding action. Then there is an affine action $\theta$ of the group $G$ on $\R^n$ and a $G$-equivariant local affine diffeomorphism 
		$$
		\text{dev}  :  M \to \R^n
		$$  
		called the {\bfseries development map} \cite{Go}. The $G$-equivariance means that for any $g\in G$ we have the following diagram
			$$
		\begin{CD}
		M  @>\text{dev}>> \R^n\\
		@VV \phi(g) V @VV \theta(g) V @.\\
		M @>\text{dev}>> \R^n
		\end{CD} 		$$
		Note that the development map can be not injective: if $\text{dev}(p)=\text{dev}(q)$ then for any $g\in G$ we have $\text{dev}\left(\phi_g (p)\right)=\text{dev}\left(\phi_g(q)\right)$.
	\end{defin}
	A flat connection $\nabla$ on a simply connected manifold $M$ sets a trivialization $TM\simeq M\times \R^n$. Since the development map is defined we can consider any tangent vector from $T(TM)=T(M\times\R^n)$ as a vector $X\times Y\in \R^n\times\R^n$. The complex structure described in the previous section acts by the rule $I(X\times Y) =-Y\times X$.

	Define an action $\psi$ of $G\ltimes_\theta \R^n$ on $TM\simeq M\times \R^n$ by 
	$$
	\psi_{a\ltimes_\theta u }(m\times v)=\phi_a(m)\times \left(\theta_a(v)+u\right),
	$$
	where $a\in G, m\in M,$ and $u, v \in \R^n,$ 
	\begin{proposition}\label{31}
		Let $G,M,\psi,\theta$ be as above. Then the complex structure $I$ on $TM$ constructed from the affine structure on $M$ is $G\ltimes_\theta \R^n$-invariant.
	\end{proposition}

	\begin{proof}
		Since there is the affine immersion $\text{dev}  :  M \to \R^n$, it is enough to prove the proposition for $M\subset\R^n$. In this case, action $\phi$ and $\theta$ define the same affine automorphism of $\R^n$ for any $g\in G$. The action $\psi$ of $G$ on any tangent vector $X\times Y\in TM\times\R^n$ is defined by 
		$$
		\psi_{g}(X\times Y)=A_g X\times A_g Y,
		$$  
		where $A_g$ is a linear part of an affine action $\phi_g$ on $\R^n$. The complex structure $I$ is defined by $I(X\times Y)=-Y\times X$. Therefore, the action $\phi$ preserves the complex structure.
	\end{proof}
	\begin{theorem}\label{T}
		Let $(M,\nabla,g)$ be a simply connected Hessian manifold, $(I,{g^\text{r}})$ the corresponding K\"ahler structure on $TM=M\times \R^n$, $G$ the group of affine isometries of $(M,\nabla,g)$, $\theta$ the action of $G$ on $\R^n$, and $\psi$ be the action of $G\ltimes_\theta\R^n$ as above. Then the K\"ahler structure $(I,{g^\text{r}})$ is invariant under the action $\psi$ of $G\ltimes_\theta\R^n$. In particular, if $\left(M,\nabla,g\right)$ is a homogeneous Hessian manifold for the group $G$ then $TM$ is a homogeneous K\"ahler manifold for the group $G\ltimes_\theta\R^n$.  
	\end{theorem}
\begin{rem}
	Theorem \ref{T} is known in the case of symmetric homogeneous Hessian manifold. Any simply connected symmetric homogeneous Hessian manifold is a direct product of flat space $\left(\R^k,\sum_{i=1}^k \left(dx^i\right)^2\right)$ and a regular convex cone $\left(V,g_{can}\right)$ (\cite{shima}). As we say in the introduction invariant K\"ahler metric on $TV\simeq V+\sqrt{-1}\R^n$ is described in \cite{VGP}.
\end{rem}

\begin{proof}[Proof of Theorem \ref{T}.]
	According to Proposition \ref{31}, the complex structure $I$ is $G\ltimes_\theta\R^n$-invariant. 
	
	By Proposition \ref{2.3}, ${g^\text{r}}$ is is locally expressed by $\text{Hess}_\mathbb{C}(4\pi^*\varphi)$. The subgroup $\{e\}\ltimes \R^n \subset G\ltimes_\theta \R^n$ acts acts trivially on the first factor of the product $M\times\R^n$. The function $\pi^*\varphi$ depends only on the first factor of $M\times\R^n$. Hence, ${g^\text{r}}=\text{Hess}_\mathbb{C}(4\pi^*\varphi)$ is invariant under the action of  $\{e\}\ltimes_\theta \R^n \subset G\ltimes_\theta \R^n$. 
	
	By \eqref{101}, 
	$
	{g^\text{r}}(X,Y)=\pi^*g(X,Y)+\pi^*g(IX,IY).
	$
	Moreover, $g$ is invariant under the action of $G=G\ltimes_\theta 0$. Hence, $\pi^*g$ is invariant under the action of $G\ltimes 0\subset G\ltimes_\theta \R^n$. Thus, ${g^\text{r}}$ is $G\ltimes 0$-invariant. Therefore, ${g^\text{r}}$ is $G\ltimes_\theta \R^n$-invariant.
	
	Actions of subgroups $\{e\}\ltimes_\theta  \R^n$ and $G\ltimes_\theta  0\subset G\ltimes_\theta \R^n$ preserve the K\"ahler structure $(I,{g^\text{r}})$. Since these subgroups generate $G\ltimes_\theta \R^n$, the K\"ahler structure $(I,{g^\text{r}})$ is $G\ltimes_\theta \R^n$-invariant.
	
\end{proof}

\section{Selfsimilar Hessian and confomally K\"ahler structures}

	\begin{defin}
		A field $\xi$ on an affine manifold $\left(M,\nabla\right)$ is called {\bfseries affine} if the flow along $\xi$ preserves the connection $\nabla$.
	\end{defin}
	\begin{proposition}\label{42}
		Let $\xi$ be an affine field on a simply connected affine manifold $(M,\nabla)$. Then there exists a vector field $\hat \xi$ on $\R^n$ such that 
		$$
		\text{dev}_* \xi = \hat \xi.
		$$
	\end{proposition}  

\begin{proof}
	The field $\xi$ defines an infinitesimal affine automorphism $A$ of $M$. Since any affine automorphisms of $M$ sets an affine automorphism $B$ of $\R^n$ (see Definition \ref{31}), $A$ sets an infinitesimal affine automorphism of $\R^n$. Since the developing map is $\text{Aut}(M,\nabla)$-equivariant, the field $\hat\xi$ of the infinitesimal affine affine automorphism $B$ satisfies $
	\text{dev}_* \xi = \hat \xi.
	$
\end{proof}
\begin{defin}
	A {\bfseries selfsimilar Hessian manifold} $\left(M,\nabla,g,\xi\right)$ is a Hessian manifold $\left(M,\nabla,g\right)$ endowed with an affine field $\xi$ satisfying $\L_\xi g =2g$. 
\end{defin}
		\begin{lemma}\label{3.4}
			Let $\left(M,\nabla,g,\xi\right)$ be a simply connected selfsimilar Hessian manifold. Consider the field $\hat\xi$ on $\R^n$ defined in Propostion \ref{42}. Denote by $\pi$ the natural projection from $M\times \R^n$ to $M$ and $\xi_1:=\xi\times0\in TM\times T\R^n$, ${\xi_2:=0\times \hat \xi\in TM\times T\R^n}$. Then we have:
		\begin{itemize}
			\item[(i)] $\L_{\xi_1}  \pi^*g=2\pi^*g$.
			
			\item[(ii)] $\L_{\xi_2} \pi^*g=0.$
			
			\item[(iii)] $\L_{\xi_1+\xi_2} I=0.$ 
		\end{itemize}
	where $\L$ is the Lie derivative. 
	\end{lemma}

	\begin{proof}
		i) Since $\pi^*g$ depends only on the projection $\pi$ on $M$, we have  
		
		$$
		\L_{\xi_1} \pi^*g = \pi^*(\L_{\pi_*\xi_1} g) =\pi^* \L_\xi g=2\pi^*g.
		$$
		
		ii) Since $\pi^* g$ depends only on the projection $\pi$ on $M$ and $\pi_*\xi_2=0$, we have
		$$
		\L_{\xi_2} \pi^*g=\pi^*\left(\L_{\pi_*\xi_2} g\right)=0.
		$$
		
		iii) Since the condition is local, we can suppose that $M\subset \R^n$. The field $\xi$ is locally defines a 1-parameter family of actions $A_t=\exp t\xi$. Then the field $\xi_1+\xi_2$ locally  defines 1-parameter family of actions $A_t\times A_t$ on $M\times \R^n\subset\R^n\times \R^n$. The action $A_t\times A_t$ preserves the complex structure. Thus $\xi_1+\xi_2=\left.\frac{d}{dt} \left(A_t\times A_t\right)\right|_{t=0}$ is holomorphic. 
	\end{proof}

\begin{defin}
	A {\bfseries selfsimilar K\"ahler manifold} $\left(M,I,g,\xi\right)$ is a K\"ahler manifold $\left(M,I,g\right)$ endowed with a holomorphic vector field $\xi$ satisfying $$\L_\xi g =2g.$$
\end{defin}
	
\begin{proposition}
	Let $(M,\nabla,g,\xi)$ be a simply connected selfsimilar Hessian manifold, $I,{g^\text{r}},\xi_1,\xi_2$ as above. Then $\left(M\times \R^n,I,{g^\text{r}}, \xi_1+\xi_2\right)$ is a selfsimilar K\"ahler manifold. 
\end{proposition}
\begin{proof}
	
	By Lemma \ref{3.4}, we have 
	$$
	\L_{\xi_1+\xi_2} \pi^*g =2\pi^*g \ \ \ \  \text{and} \ \  \ \ \L_{\xi_1+\xi_2} I=0.
	$$
	Combining this with \eqref{101}, we obtain that
	$$
	\L_{\xi_1+\xi_2} {g^\text{r}}= 2 {g^\text{r}}.
	$$
	By the item (iii) of Lemma \ref{3.4}, $\xi_1+\xi_2$ is a holomorhic vector field. Thus,  ${(M\times\R^n,I,{g^\text{r}}, \xi_1+\xi_2)}$ is a selfsimilar K\"ahler manifold. 
\end{proof}

\begin{lemma}\label{36}
	Let $(M,\nabla,g,\xi)$ be a simply connected selfsimilar Hessian manifold, $(I,{g^\text{r}})$ the corresponding K\"ahler structure on $TM=M\times\R^n$, and $G$ the group of affine isometries of $M$ preserving $\xi$. Then $\pi^*g(\xi_1,\xi_1)$ is a $G\ltimes_\theta \R^n$-invariant function on $M\times\R^n$ satisfying $$\L_{\xi_1+\xi_2}(\pi^*g(\xi_1,\xi_1))=2\pi^*g(\xi_1,\xi_1)$$ at any point.
\end{lemma}

\begin{proof}
	Since $\pi^*g$ and $\xi_1$ are $G\ltimes_\theta \R^n$-invariant, the function $\pi^*g(\xi_1,\xi_1)$ is $G\ltimes \R^n$-invariant. A Lie derivative of a metric if defined by
	$$
	\left(\L_X g\right)(X,Y)=\left(\L_X\right)g(X,Y)-g\left(\L_X Y, Z\right)-g\left(Y,\L_X Z\right).
	$$  
	We have $[\xi_1,\xi_2]=0$. Therefore, 
	$$
	\L_{\xi_1+\xi_2}\left(\pi^*g(\xi_1,\xi_1)\right)=\left(\L_{\xi_1+\xi_2} \pi^*g\right)\pi^*g(\xi_1,\xi_1).
	$$
	Combining this with items (i) and (ii) of Lemma \ref{3.4}, we get
	$$
	\L_{\xi_1+\xi_2}\left(\pi^*g(\xi_1,\xi_1)\right)=2 \pi^*g(\xi_1,\xi_1).
	$$

\end{proof}
\begin{defin}
	A {\bfseries (globally) conformally K\"ahler structure} on a manifold $N$ is a collection $\left(I,f\omega\right)$, where $\left(I,\omega\right)$ is a K\"ahler structure and $f$ is a positive definite function on $N$.
\end{defin}

	Let $G$ be the group of affine isometries of a selfsimilar affine manifold $\left(M,\nabla,\xi,g\right)$ preserving $\xi$ and $\xi$ is complete. Let $\rho_t=\left(\exp t\xi\right)$ be an affine action of $\R$ on $M$. The action of $G$ preservers $\xi$. Therefore, this action commutes with $\rho$. Thus, there is an action of $\R\times G$ on $M$ defined by $\tilde\theta_{t\times g}(m)=\rho_t\left(\theta_g (m)\right)$.
	
\begin{theorem}\label{46}
	Let $(M,\nabla,g,\xi)$ be a simply connected selfsimilar Hessian manifold, $(I,\omega)$ the corresponding K\"ahler structure on $M\times \R^n$, $G$ the group of affine isometries of $(M,\nabla,g)$ preserving $\xi$, the actions $\theta$ and $\tilde \theta$ as above, $\pi :M\times \R^n\to M$ the natural projection and
	${\omega_{c.K.}=\pi^*\left(g(\xi,\xi)^{-1}\right)\omega}$. Then the following conditions are satisfies: 
	\begin{itemize}
		\item [i)] The pair $\left(I,\omega_{c.K.} \right)$ is a $G\ltimes_\theta\R^n$-invariant conformally K\"ahler structure on $M\times\R^n$.
		\item [ii)] If the field $\xi$ is complete then $\left(I,\omega_{c.K.} \right)$ is $\left(\R \times G\right)\ltimes_{\hat \theta} \R^n$-invariant.
		\item [iii)] If the field $\xi$ is complete and  $G$ acts transitively on the level surface $S=\{g(\xi,\xi)=1\} \subset M$ then $\left(M\times \R^n,I,\omega_{c.K.} \right)$ is homogeneous conformally K\"ahler manifold for the group $\left(\R \times G\right)\ltimes_{\tilde \theta} \R^n$.
	\end{itemize}
\end{theorem}

\begin{proof}
	i) According to Lemma \ref{36} the function $\left(\pi^*g(\xi_1,\xi_1)\right)^{-1}$ is $G\ltimes \R^n$-invariant. According to Theorem \ref{T}, $\omega$ is $G\ltimes \R^n$-invariant. Hence, $\omega_c.K$ is  $G\ltimes \R^n$-invariant. By Proposition $\ref{31}$, the complex structure $I$ is $G\ltimes \R^n$-invariant. 
	Thus, the conformally K\"ahler structure $\left(I,\omega_{c.K.}\right)$ is $G\ltimes \R^n$-invariant.

	ii) The action of a subgroup
	$$
	\left(\R \times 0\right)\ltimes_{\hat \theta} 0\subset \left(\R \times G\right)\ltimes_{\hat \theta} \R^n.
	$$
	is generated by the vector field $\xi_1+\xi_2$. By Lemma \ref{36}, we have
	$$
	\L _{\xi_1+\xi_2} \left(\pi^*g(\xi_1,\xi_1)\right)^{-1}=\left(\pi^*g(\xi_1,\xi_1)\right)^{-2} \L_{\xi_1+\xi_2} \left(\pi^*g(\xi_1,\xi_1)\right)^{-1}=-2 \left(\pi^*g(\xi_1,\xi_1)\right)^{-1}.
	$$ 
	Therefore,
	$$
	\L_{\xi_1+\xi_2}\omega_{c.K.} =\L_{\xi_1+\xi_2} \left(\pi^*\left(g(\xi,\xi)^{-1}\right)\omega\right)=\L_{\xi_1+\xi_2}\left(\pi^*g(\xi_1,\xi_1)\right)^{-1}\omega+\L_{\xi_1+\xi_2}\omega=-2\omega+2\omega=0.
	$$
	Hence, if the field $\xi$ is complete then $\omega_{c.K.}$ is invariant under the action of the subgroup 
	$
	{\left(\R \times 0\right)\ltimes_{\hat \theta} 0}.
	$ 
	By Lemma \ref{3.4}, the complex structure $I$ is invariant under the action of ${\left(\R \times 0\right)\ltimes_{\hat \theta} 0}$. As we showed above, the c.K. structure $\left(I,\omega_{c.K.}\right)$ is invariant under the action of the group $${G\ltimes_\theta \R^n=\left(0 \times G\right)\ltimes_{\hat \theta} \R^n}.$$
	Therefore, $\left(I,\omega_{c.K.}\right)$ is $\left(\R \times G\right)\ltimes_{\hat \theta} \R^n$-invariant.
	
	iii) The subgroups $\left(\R\times \{e\}\right)\ltimes_{\tilde{\theta}} \{0\}\subset \left(\R \times G\right)\ltimes_{\hat \theta} \R^n$ act on the function $\pi^* g(\xi,\xi)$ be the rule
	$$
	\psi^*_{\left(r\times e\right)\ltimes_{\tilde{\theta}} 0} \pi^* g(\xi,\xi)=r^2 \pi^* g(\xi,\xi).
	$$
	Thus, we can can send any point of $M\times\R^n$ by the action of this group to a point of $S\times\R^n$. The group $\left(\{0\} \times G\right)\ltimes_{\tilde \theta} \R^n=G\ltimes_\theta\R^n$ acts transitively on $S\times \R^n$. Therefore, $\left(\R \times G\right)\ltimes_{\tilde \theta} \R^n$ acts transitively on $M\times\R^n$.
\end{proof}
	\section{An example: homogeneous regular convex cones}
	\begin{defin}
		A set $V\in \R^n$ is a {\bfseries cone} if and only if for any $x\in V$ and $a\in \R^{>0}$ we have $ax\in V$. We say that a cone $V$ is {\bfseries regular} if it does not contain any straight full line. 
	\end{defin}

	\begin{theorem}[\cite{vinb},\cite{VGP}]\label{52}
		Let $V$ be a homogeneous convex cone. Then there exists a function $\varphi$ satisfying the following conditions.
		\begin{itemize}
			\item [(i)] The bilinear form $g_{can}=\text{Hess} (\ln \varphi)$ is $\text{Aut(V)}$-invariant Hessian metric.
			\item [(ii)] The level line $\varphi=1$ is $\text{Aut}_{\text{SL}}(V)$-invariant.
			\item [(iii)] The bilinear form ${g^\text{r}}_{can}$ is $\text{Aut}(V)\ltimes \R^n$-invariant K\"ahler structure on $TM\simeq V\oplus \sqrt{-1}\R^n$. 
			\item [(iv)] The bilinear form $g_{con} =\text{Hess}\  \varphi$ is $\text{Aut}_{\text{SL}}(V)$-invariant Hessian metric satisfying the condition ${\L_\rho g_{con} =-ng}$, where $\rho=\sum_{i=1}^n x^i\dxi$ is the radiant vector field.
		\end{itemize}
	\end{theorem}

According to item (iv) of Theorem \ref{52}, $\left(V,g_{con},-\frac{2}{n}\rho\right)$ is a selfsimilar Hessian manifold. 
	\begin{defin}
		Let $V$ be a convex regular cone. Then $V\oplus \sqrt{-1} \R^n$ is biholomorphic to a bounded domain in $\CC^n$. All complex domains arising by this way are called {\bfseries Siegel domain of the first kind}.
	\end{defin}
 	The item $(iii)$ of Theorem \ref{52} states that any homogeneous Siegel domain of the first kind admits a homogeneous K\"ahler structure. Applying Theorem \ref{46} to the manifold $V$ with a metric $g_{con}$ and action of $\text{Aut}_{\text{SL}}(V)$ we get a similar result.
 	\begin{cor}
 		Any homogeneous Siegel domain of the first kind admits a homogeneous conformally K\"ahler structure.
 	\end{cor}

\section{The c-map for homogeneous manifolds}
\begin{defin}
	A  {\bfseries special K\"ahler manifold} $\left(M,I,g,\nabla\right)$ is a K\"ahler manifold  $\left(M,I,g\right)$ endowed with a torsion free symplectic connection $\nabla$ such that such that $g$ is a Hessian metric with respect to $g$. 
\end{defin}
\begin{theorem}[\cite{CMMS2}, \cite{ACM}]
Let $\left(M,I,g,\nabla\right)$ be a special K\"ahler manifold. Using the connection $\nabla$ we can identify $T\left(T^*M\right)=T^h\left(T^*M\right)\oplus T^v\left(T^*M\right)=\pi^*TM\oplus \pi^* TM$, where $\pi: TM\to M$ is the canonical projection, $T^h\left(T^*M\right)$ and $T^v\left(T^*M\right)$ are the horizontal and vertical distribution defined by $\nabla$. Using this identification we define 
$$
g^\text{c}=\begin{pmatrix}
g & 0 \\
0 & g^{-1}
\end{pmatrix}, \ \ \  \
I_1=
\begin{pmatrix}
I & 0 \\
0 & I^*
\end{pmatrix}, \ \ \ \
I_2=
\begin{pmatrix}
0 & -\omega^{-1} \\
\omega & 0
\end{pmatrix}, \ \ \ \
I_3=I_1 I_2.
$$	
Then $\left(T^*M,I_1,I_2,I_3, {g^\text{c}}\right)$ is a hyper-K\"ahler manifold.
\end{theorem}  
\begin{defin}
	The correspondence which associates the hyper-Kähler manifold $\left(T^*M,I_1,I_2,I_3, {g^\text{c}}\right)$ to a special K\"ahler manifold $\left(M,I,g,\nabla\right)$ is called the {\bfseries (affine) c-map} (see \cite{AC}).
\end{defin}
The flat connection $\nabla$ sets a trivialization $TM\simeq M\times\R^n$. An affine action $\phi$ of $G$ on $M$ sets an affine action of $G$ on $\R^n$ (see section 3). Since the bilinear form $\omega$ is $\nabla$-flat, we can consider $\omega$ as a bilinear form on $\R^n$. Hence, $\omega$ sets an identification of $\R^n$ and $\left(\R^n\right)^*$. Thus, we can consider the conjugate action $\hat\theta$ of $G$ on $\left(\R^n\right)^*$.
Define the action $\hat\psi$ of $G\ltimes_{\hat\theta} \left(\R^n\right)^*$ on $T^*M\simeq M\times \left(\R^n\right)^*$ by 
$$
\hat\psi_{a\ltimes_\theta u }(m\times v)=\phi_a(m)\times \left(\hat\theta_a(v)+u\right),
$$
The $\nabla$-flat bilinear form $\omega$ sets an identification $TM\simeq M\times \R^n\simeq M\times \left(\R^n\right)^*\simeq T^*M$. Thus, we can define an action of $G\ltimes_{\theta^*} \left(\R^n\right)^*=G\ltimes_\theta \R^n$ on $T^*M\simeq M\times\left(\R^n\right)^*$.
and fields $\xi_1\in T^h\left(T^*M\right)$, $\xi_2\in T^v\left(T^*M\right)$. 

\begin{theorem}\label{64}
	Let $\left(M,I,g,\nabla\right)$ be a special K\"ahler manifold, $G$ the group of affine holomorphic isometries of $\left(M,I,g,\nabla\right)$, $\hat \theta$ the action of $G$ on $\left(\R^n\right)^*$, $\hat\psi$ the action of $G\ltimes_\theta \left(\R^n\right)^*$ on $ \left(\R^n\right)^*$ as above. Then the hyper-K\"ahler structure  $\left(I_1,I_2,I_3, {g^\text{c}}\right)$ is invariant under the action $\hat\psi$ of $G\ltimes_{\hat \theta}\R^n$. In particular, if $\left(M,I,g,\nabla\right)$ is a homogeneous special K\"ahler manifold for a group $G$ then $\left(T^*M, I_1,I_2,I_3, {g^\text{c}}\right)$ is a homogeneous hyper-K\"ahler manifold for the group  $G\ltimes_{\hat\theta} \left(\R^n\right)^*$.
\end{theorem}

\begin{proof}
	Tensors $g^c, I_1,I_2,I_3$ depends only on the first factor of $M\times\left(\R^n\right)^*$. The action $\hat\psi$ acts on the first factor as $\phi$. Tensors $g,I,\omega$ on $M$ are invariant with respect to the action $\phi$ of the group $G$. Hence, $g^{-1}, I^*,-\omega^{-1}$ are invariant with respect to the action $\phi$, too. Therefore, $g^c,I_1,I_2,I_3$ are invariant with respect to the action $\hat\psi$. Moreover, if the action $\phi$ of $G$ on $M$ is transitive, then the action  $\hat\psi$  of $G\ltimes_\theta \left(\R^n\right)^*$ on $ T^*M\simeq M\times\left(\R^n\right)^*$ is transitive. Thus, if $\left(M,I,g,\nabla\right)$ is a homogeneous special K\"ahler manifold for a group $G$ then $\left(T^*M, I_1,I_2,I_3, {g^\text{c}}\right)$ is a homogeneous hyper-K\"ahler manifold for the group  $G\ltimes_{\hat\theta} \left(\R^n\right)^*$.
\end{proof}

\section{Selfsimilar special K\"ahler and conformally hyper-K\"ahler manifolds}

\begin{defin}
		A {\bfseries (globally) conformally hyper-K\"ahler structure} on a manifold $N$ is a collection $\left(f g, I_1,I_2,I_3\right)$, where $\left(g, I_1,I_2,I_3\right)$ is a hyper-K\"ahler structure and $f$ is a positive definite function on $M$. 
\end{defin}

\begin{defin}
	A {\bfseries selfsimilar special K\"ahler manifold} $\left(M,I,g,\nabla,\xi\right)$ is a special K\"ahler manifold $\left(M,I,g,\nabla\right)$  endowed with an affine vector field $\xi$ satisfying $\L_\xi g= 2g$.
\end{defin}
\begin{theorem}\label{73}
	Let $\left(M,I,g,\nabla,\xi\right)$ be a selfsimilar special K\"ahler manifold,  $G$ the group of affine holomorphic isometries of  $\left(M,I,g,\nabla\right)$ preserving $\xi$, $\left(I_1,I_2,I_3, {g^\text{c}}\right)$ the corresponding hyper-K\"ahler structure on 
	$T^*M$, $\pi :M\times \R^n\to M$ the natural projection, and $g_{c.h.K}=\pi^*\left(g(\xi,\xi)\right)^{-1}g^c$. Then the following conditions are satisfies:
	
		\begin{itemize}
		\item [i)] The collection $\left(g_{c.h.K},I_1,I_2,I_3\right)$ is 
		$G\ltimes_{\hat\theta}\left(\R^n\right)^*$-invariant hyper-K\"ahler structure on $T^*M$.
		\item [ii)] If the field $\xi$ is complete then $\left(g_{c.h.K},I_1,I_2,I_3\right)$ is $\left(\R \times G\right)\ltimes_{\hat \theta} \left(\R^n\right)^*$-invariant.
		\item [iii)] If the field $\xi$ is complete and  $G$ acts transitively on the level surface $S=\{g(\xi,\xi)=1\} \subset M$ 
		then $\left(g_{c.h.K},I_1,I_2,I_3\right)$ is homogeneous conformally hyper-K\"ahler manifold for the group 
		$\left(\R \times G\right)\ltimes_{\tilde \theta} \left(\R^n\right)^*$.
	\end{itemize}
\end{theorem}
\begin{proof}
	By the same argument as in the Theorem \ref{64}, the hyper-K\"ahler structure $\left(I_1,I_2,I_3, {g^\text{c}}\right)$ is $G\ltimes_{\hat\theta}\left(\R^n\right)^*$-invariant. Then the proof of Theorem\ref{73} coincides with the poof of Theorem \ref{46}.

%	Since, the $g$ metric on $M$ and the field $\xi$ are $G$-invariant and the subgroup $\{e\}\ltimes \left(\R^n\right)^*\subset G\ltimes \left(\R^n\right)^*$, preserves the fibers of the cotangent bundle $T^*M\simeq M\times \left(\R^n\right)^*$  over $M$, the bilinear form $\pi^*g$ and the field $\pi^*\xi\in T^h\left(T^*M\right)=\pi^*TM$ are $G\ltimes_{\hat \theta} \left(\R^n\right)^*$-invariant. Therefore, the function $\pi^*\left(g(\xi,\xi)\right)=\left(\pi^*g\right)\left(\pi^*\xi,\pi^*\xi\right)$ is $G\ltimes_{\hat \theta} \left(\R^n\right)^*$-invariant. Thus, $\left(g_{c.h.K},I_1,I_2,I_3\right)$ is 
	%$G\ltimes_{\hat\theta}\left(\R^n\right)^*$-invariant conformally hyper-K\"ahler structure.

\end{proof}

\end{document}